\documentclass[12pt,draft]{article}

\usepackage{times}
\usepackage{a4,amssymb}
\def\Bbb{\mathbb}

\title{\bf Higher cotangent cohomology of rational surface singularities}
\author{Jan Stevens\thanks{Partially supported by the 
               Swedish Research Council (Vetenskapsr\aa det).}}
\date{}

\def\m{{\mathfrak m}}
\def\wtx{\widehat{X}}
\def\wt{\widehat}
\def\wl#1{\overline{#1}}
\def\sier#1{{\cal O}_{#1}}
\def\C{{\Bbb C}}  \def\L{{\Bbb L}}  \def\Z{{\Bbb Z}} \def\Q{{\Bbb Q}}
\def\cA{{\cal A}}
\def\cM{{\cal M}} \def\cT{{\cal T}} \def\cN{{\cal N}} \def\cF{{\cal F}}
\def\cP{{\cal P}} \def\cQ{{\cal Q}} \def\cR{{\cal R}} \def\cS{{\cal S}}

\def\al{\alpha}
\def\kn{i}
\def\vp{\varphi}
\def\phi{\varphi}
\def\kb{{\scriptscriptstyle \bullet}}
\def\Rank{\mathop{\rm Rank}}
\def\Ker{\mathop{\rm Ker}}
\def\Coker{\mathop{\rm Coker}}
\def\cg{\mathop{\rm cg}}
\def\cod{\mathop{\rm cod_{AC}}}
\newcommand{\gHom}{\mbox{\rm Hom}}
\newcommand{\Harr}{\mathop{\rm Harr}\nolimits}
\def\lra{\longrightarrow}
\def\mapright#1{\mathrel{\mathop{\longrightarrow}\limits^{#1}}}

\def\roep #1.{\medbreak\noindent{\sl #1\/}.\enspace}
\def\endroep{\par\medbreak}
\def\thesection{\arabic{section}.}
\def\thesubsection{\thesection\arabic{subsection}}
\long\def\proclaim#1. #2\par{\medbreak\noindent{\bf #1.}\enspace{\sl#2}\par
\ifdim\lastskip<\medskipamount \removelastskip\penalty55\medskip\fi}
\def\Box{\square}
\newcommand{\qed}{{
\unskip\nobreak\hfil\penalty50\hskip2em\hbox{}\nobreak\hfil$\Box$
\parfillskip=0pt \finalhyphendemerits=0 \par\bigbreak}}
\def\mref#1{{\rm(\ref{#1})}}

\newcommand{\neu}[1]{\protect\refstepcounter{subsection}\protect
   \label{#1}\vspace{1ex}
   {\bf (\thesubsection)} \enspace\ignorespaces}

\parskip=\medskipamount

\begin{document}

\maketitle

The modules $T^1$ and $T^2$ play an important role in deformation
theory, the first as space of infinitesimal deformations, while the
obstructions land in the second. 
Much work has been done to compute their dimension for rational
surface singularities, culminating in the formulas of 
{\sc Christophersen} and {\sc Gustavsen} \cite{CG}.
The `correct' way to define $T^1$ and $T^2$ also yields higher
$T^i$. The purpose of this note is to generalise the
dimension formulas to these modules. 

For rational surfaces with reduced fundamental cycle 
{\sc De Jong} and {\sc Van Straten} \cite{dJvS}
showed that the dimension of $T^2$ can be computed inductively from the
multiplicities of the singularities on the successive blow-ups.
In the general case of a rational surface singularity $X$ of multiplicity 
$d\geq 3$
with first blow-up $\wtx$ the result \cite[Thm.~3.8]{CG} is: 
$$
\dim T^2_X= (d-1)(d-3) + \dim T^2_{\wtx} +c(X)\;,
$$
where $c(X)$ is an unavoidable correction term,
which vanishes for large classes of singularities (including 
reduced fundamental cycle). With the  same correction term one has
\cite[Thm.~3.13]{CG}:
$$
\dim T^1_X= (d-3) + \dim T^1_{\wtx} +c(X)\;.
$$
Our main result is 
that the corresponding formula holds for the higher $T^i$
without correction term:
$$
\dim T^i_X= f_i(d) + \dim T^i_{\wtx}
$$
with $f_i(d)$ an explicit function of $d$. 
If $\wtx$ is smooth, then $\dim T^i_X$ equals $f_i(d)$,
so the value of $f_i(d)$ follows from the computations in \cite{AS}
for the cone over the rational normal curve
of degree $d$; see also \mref{prel-rnc}.

We follow the arguments of \cite{CG}  closely, replacing the
computations of the different $T^2$'s in terms of functions
on relations by computations with Harrison cohomology.
By means of a Noether normalisation (a flat map $X \to S$
of degree $d$ onto a smooth surface) it suffices to look at
relative Harrison cohomology, whose defining complex is
$\sier S$-linear and therefore much smaller. We use here
in an essential way that the iterated hyperplane section
of a rational surface singularity is the fat point of 
minimal multiplicity
$d$. 
The computations with relative Harrison cohomology work well
for $i\geq 2$. 
Our arguments reprove the $T^2$-formula from \cite{CG}.

{\sc De Jong} and {\sc Van Straten} \cite{dJvS} prove their
formulas for $T^1$ and $T^2$  with
a one-parameter deformation
of the singularity (with reduced
fundamental cycle) to the cone over the rational normal curve
of degree $d$ and all singularities on the first blow up. 
In this deformation the codimension of the Artin component
and the dimension of $T^2$ 
are  constant. By induction on the singularities on the blow up
it then follows that the number of equations
needed for the base space equals the dimension of $T^2$, a property known
as the surjectivity of the obstruction map.
Surjectivity still follows if we forget about rational double points and 
relax the requirements. I had hoped to prove that every rational surface
singularity admits such a {\sl good maximal deformation\/}, but examples
show that this is not the case. 

I am grateful to Jan Christophersen for discussions and comments.

\section{Cotangent cohomology}
\neu{prel-gen}
The definition of $T^i_X$ is most conveniently given in the context of
a more general theory and we set therefore $T^i_X:=T^i(X/{\rm pt};\sier X)$.
Let $X\to S$ be a map of analytic germs  
and $M$ an $\sier X$-module.
One gets the  {\sl cotangent cohomology groups\/} as
$$
T^i(X/S;M):=H^i\!\left(\gHom_{\sier X}(\L^{X/S}_*,M)\right)
$$
with $\L^{X/S}_*$ being the so-called cotangent complex. What we need
to know about it is summarised in \cite{BC}, \cite{CG}. 
A good reference, also for the next subsection, is
\cite{L}. 
We actually work with the analytic version of the cotangent complex,
which can be constructed from a Tyurina resolvent of the analytic
$\sier S$-algebra $\sier X$,  
see {\sc Palamodov}'s 
survey \cite{P3} and for more
details his papers \cite{P4}, \cite{P2} and \cite{P}.

\neu{prel-harr}
The higher cotangent cohomology can also be computed as 
Harrison cohomo\-logy. To give the definition 
we first recall Hochschild cohomology.
Let $A$ be a commutative algebra of essentially finite type over a base
{\sl ring\/} $k$.
For an $A$-module $M$, the
{\sl Hochschild cohomology} $HH^i(A/k;M)$ is the cohomology of the complex
$$
C^i(A/k;M):=\gHom_k(A^{\otimes i}, M)
$$
with differential 
$$
\displaylines{\;\;
(\delta f) (a_0,\dots,a_i):={}\hfill\cr\hfill
 a_0f(a_1,\dots,a_i)+
\sum_{j=1}^i(-1)^j f(a_0,\dots,a_{j-1}a_{j},\dots,a_i)
+(-1)^{i+1}a_i f(a_0,\dots,a_{i-1})\;.
}
$$ 
The same cohomology can also be obtained 
from the so-called {\it reduced\/} subcomplex $\wl C^\kb(A/k;M)$
consisting only of those maps $f\colon A^{\otimes i}\to M$ that vanish 
whenever at least one of the arguments equals 1.

\roep Definition.
A permutation $\sigma$ is called a {\sl $(p,q)$-shuffle\/} 
if $\sigma(1)<\dots<\sigma(p)$ and
$\sigma(p+1)<\dots<\sigma(p+q)$. Moreover, in the group algebra
$\Z[S_{p+q}]$ we define the elements
$$
\mbox{\rm sh}_{p,q}:=\sum_{(p,q){\rm-shuffles}} 
\mbox{sgn}(\sigma)\, \sigma\;.
$$

These elements give rise to the 
so-called shuffle invariant subcomplexes
$$
C_{\rm sh}^i(A/k;M):=
\left\{f\in \gHom_k(A^{\otimes i}, M)\mid
f\circ{\rm sh}_{p, i-p}=0\mbox{ for  }0<p<i\right\}
$$
and $\wl{C}_{\rm sh}^i(A/k;M)\subset\wl{C}^i(A/k;M)$
defined in the same manner. Both complexes yield the same cohomology, which is 
called  {\sl Harrison cohomology\/} 
$\Harr^i(A/k;M)$.

\proclaim Proposition.
If\/ $\Q\subset k$, then 
Harrison cohomology 
is a direct summand of  Hoch\-schild cohomology.

In the analytic case one has to use the analytic tensor product; precise
definitions are given by Palamodov \cite{P}. We can now
state the relation with 
cotangent cohomology.

\proclaim Theorem.
If the map $X\to S$ is  {\sl flat}, then
$$
T^i(X/S;M) \cong \Harr^{i+1}(\sier X/\sier S;M)\;.
$$

\neu{prel-norm}
Let $X$ be a Cohen-Macaulay singularity of dimension $n$ and multiplicity
$d$. 
We can choose as  {\sl Noether normalisation\/} $X\to \C^n$ a 
flat map of degree $d$. In fact, this is only a map of germs of 
analytic spaces, but we are abusing notation here.
The exact sequence relating absolute and relative cotangent
cohomology gives for an $\sier X$-module $M$
$$
\def\longrightarrow{\to}
\cdots \longrightarrow T^i(X/\C^n;M) \longrightarrow T^i(X;M)
\longrightarrow T^i(\C^n;M)\longrightarrow T^{i+1}(X/\C^n;M)
\longrightarrow \cdots
$$
As $\C^n$ is smooth, $T^i(\C^n;M)=0$ for $i\geq1$.
We obtain therefore
\proclaim Lemma.
For  any $\sier X$-module $M$ one has 
$T^i(X;M) \cong T^i(X/\C^n;M)$ for $i\geq 2$.

In particular, if $X$ has minimal multiplicity $d$ (e.g., 
a rational surface singularity),
so  $\mbox{\rm embdim}\,X=d+n-1$, we may choose coordinates 
$(z_1,\dots,z_{d+n-1})$  such that the 
projection on the space spanned by the last $n$ coordinates
is a Noether normalisation.

In terms of rings we have a regular local ring $P=\C\{z_d,\dots,z_{d+n-1}\}$
and a homomorphism $P\to A$ making the local ring $A$ of $X$ into
a free $P$-module of rank $d$, with basis $\{1,z_1,\dots,z_{d-1}\}$. 
The equations for $X$ are of the form
$$
z_iz_j= \sum_{\nu=1}^{d-1}p_\nu z_\nu+ p_0, \qquad 1\leq i,j\leq d-1\;,
$$
where $p_\nu\in P$ for $\nu=0$, \dots, $d-1$.

\neu{prel-zero}
Let $X$ be a rational surface singularity of multiplicity $d\geq 3$
with local ring $A$ and Noether normalisation $X\to \C^2$, or in terms
of rings $P\to A$.
We can obtain the higher cotangent cohomology as
Harrison cohomology $\Harr^{i+1}(A/P;M)$ computed
from the reduced complex. 
As illustration of this technique we repeat here the
proof of the following lemma, shown in \cite{AS}.
We shall need the argument later on.

\proclaim Lemma.
The natural map $\Harr^{i+1}(A/P;A)\to \Harr^{i+1}(A/P;\C)$ is the zero map.

\roep Proof.
As $\{1,z_1,\dots,z_{d-1}\}$ forms a basis of $A$ as $P$-module, 
a reduced
Harrison $(\kn+1)$-cocycle $f$ is, by $P$-linearity, determined
by its values on the $(\kn+1)$-tuples of the coordinates $z_1$, \dots,
$z_{d-1}$.
Suppose $f(z_{j_0}, \dots, z_{j_\kn}) \notin \m_A$. 
Since $d\geq 3$, we may choose a $z_k$ with
$k\in\{1,\dots,d-1\}$ and $k\neq j_0$. Hence,
$$
\displaylines{
\qquad
0=(\delta f)(z_{j_0}, \dots, z_{j_\kn},z_k)=
z_{j_0}f(z_{j_1}, \dots, z_k) \pm f(z_{j_0}, \dots, z_{j_\kn})z_k
\hfill\cr\hfill{}+
\mbox{ terms containing products $z_iz_j$ as arguments}\;.
\qquad
}
$$
Using the equations for $X$ as in \mref{prel-norm} we may
again apply $P$-linearity to see that the latter terms are contained in
$\m_P\cdot A$. Hence, modulo $\m_P$, these terms vanish, but
the resulting equation inside 
$V=\m_A/\m_PA$ contradicts the
fact that $z_{j_0}$ and $z_k$ are linearly independent.
\qed

\neu{prel-hps}
Another powerful tool is the comparison between the cotangent cohomology
of a singularity and its hyperplane section \cite{BC}.
Let $f\colon X\to\C$ be a flat map
such that both $X$ and the special fibre $H$ have isolated singularities.

\proclaim Main Lemma {\rm(\cite{BC}, (1.3.2))}.
There is a long exact sequence
$$
T^1(X/\C;\sier X) \lra
T^1_H \longrightarrow
T^2_X \mapright f
T^2_X \longrightarrow
T^2_H \longrightarrow
T^3_X \mapright f
T^3_X \longrightarrow
\cdots\;.
$$
Moreover, $\dim\,\raisebox{0.4ex}{$T^2_X$}\left/
            \raisebox{-0.4ex}{$f\cdot T^2_X$}\right.= \tau_H-e_f$
with $\tau_H:=\dim\,T^1_H$ and $e_f$  the dimension of the smoothing
on which $f$ lies in the versal base space of $H$.

\neu{prel-rnc}
The Main Lemma becomes particularly useful if the maximal ideal of
$X$ annihilates the $T^i_X$ for $i\geq 2$. In \cite{AS} it is shown
that the cone over the rational normal curve of degree $d$ 
has this property. Moreover a formula is given for the dimension
of all $T^i$.

To formulate it we first consider the fat point $Z_m$  
of minimal multiplicity $d=m+1$, which is 
the iterated hyperplane section of the cone over the rational normal
curve. We define the number $c_{m,k}$ as the dimension
of $T^{k-1}({Z_m};\C)$.
The dimension of  $T^i_{Z_m}$ is then $mc_{m,i+1}-c_{m,i}$.
By e.g.~\cite[2.4]{AS}
$$
c_{m,k}= \frac 1k \sum_{q|k}(-1)^{k+\frac kq}\mu(q)m^{\frac kq}\;.
$$
The first few values are
$c_{m,1}  =   m $,
$c_{m,2}  =  (m^2+m)/2   $,
$c_{m,3}  =  (m^3-m)/3 $,
$c_{m,4}  =  (m^4-m^2)/4  $,
$c_{m,5}  =  (m^5-m)/5 $,
$c_{m,6}  =  (m^6+m^3-m^2-m)/6$.
We combine them in the power series
$$
Q_d(t)=  \sum_{i\geq1} c_{d-1,i}t^i\;.
$$

\proclaim Proposition {\rm\cite[4.7]{AS}}.
The Poincar\'e series
$$
P_d(t)=\sum_{i\geq 1}f_i(d)\cdot t^i
:=\sum_{i\geq 1} \dim\, T^i_{X_d}\cdot t^i
$$
of the cone $X_d$ over the rational normal curve is given by
$$
P_d(t)\,=\;\left(Q_d(t)+2t+2\right)\cdot
\frac{(d-1)t-t^2}{(t+1)^2}
\;-\; \frac{2t}{t+1}\,.
$$

For low values of $i$ we get as dimensions of $T^i_{X_d}$:
%
%
\begin{eqnarray*}
f_1(d) &=& 2d-4 \\
f_2(d) &=& (d-1)(d-3)   \\                          
f_3(d) &=& (d-1)(d-2)(d-3)/2  \\              
f_4(d) &=& (d-1)(d-2)(2d^2-8d+9)/6 \\
f_5(d) &=& (d-1)(d-2)^2(3d^2-8d+9)/12 \\
f_6(d) &=& (d-1)(d-2)(12d^4-66d^3+153d^2-179d+90)/60\;.
\end{eqnarray*}
For the cone $X_d$ the dimension of the $T^i$ equals the number 
of generators as $\sier X$-module. For an $\sier X$-module $M$ we denote the 
minimal number of generators, which is $\dim M/\m M$, by $\cg M$.
Feeding the above results in the  Main Lemma and using the arguments of 
\cite[5.1]{BC}, one obtains (see \cite{AS}):

\proclaim Theorem.
For all rational surface singularities of multiplicity $d$ one has
that $\cg T^i=f_i(d)$ for $i\geq 2$.


\section{The tangent complex with values in a sheaf}
\neu{blow-cg}
To describe the relation between cotangent cohomology of
a singularity and its first blow-up one has to globalise 
local constructions. 
We use the analytic cotangent complex of {\sc Palamodov};
for a general overview see \cite{P3}, while some technical details
are to be found in \cite{P2}.

Let $f\colon X\to S$ be a map of complex spaces. An analytic sheaf over $X$
is  a morphism of complex spaces $\pi\colon Y\to X$ together with an
$\sier Y$-sheaf $\cF$. 
On $Y$ we have the sheaf  $\cT^i(X/S;\cF)$ whose stalk at a point 
$y\in Y$ is the cotangent cohomology $T^i(\sier {X,\pi(y)}/ \sier {S,f\pi(y)};
\cF_y)$. This sheaf occurs in a local to global
spectral sequence
$$
E_2^{pq}=R^p(f\pi)_*\left(\cT^q(X/S;\cF)\right) \Longrightarrow 
T^i(X/S;\cF)\;.
$$

In our application, where we will calculate the left hand side,  also the
right hand side reduces to a familiar object. Under the condition
that $R^q f_*(\cF)=0$ for $q>0$ {\sc Christophersen} and {\sc Gustavsen}
\cite{CG}
obtain from Prop.~56 of the Appendix in \cite{An} that 
$T^i(X/S;\cF) \simeq T^i(X/S;f_*\cF)$, where $f_*\cF$ is a sheaf on $X$.
For lack of reference we give here a proof using 
{\sc Palamodov}'s theory in the case that $f\colon X\to S$
is a finite map of germs.

\neu{cover}
A {\sl polyhedron\/} $P$ in a complex space $X$ is a relatively compact subset
$P\subset U \subset X$ together with a proper embedding
$\phi \colon U \to V \subset \C^N$, where $V$ is an open neighbourhood of
the unit polydisc $D_N$, such that $P=\phi^{-1}(D_N)$.
A polyhedral covering $\cP=\{P_\al,\al\in \cA\}$ should satisfy
$\wl P_\beta \subset U_\al$, if $\wl P_\beta \cap \wl P_\al \neq
\emptyset$. A simplex of the {\sl nerve\/} $\cN(\cP)$ is 
a map $A\colon [n] \to \cA$, written $A=(\al_0,\dots,\al_n)$, such
that $\bigcap \wl P_{\al_i} \neq\emptyset$. The set $P_A= \bigcap  P_{\al_i}$
is again a polyhedron with $\phi_A=\prod \phi_{\al_i} \colon
\bigcap U_{\al_i} \to \prod V_{\al_i}\subset \C^{N_A}$, where $N_A=
\sum N_{\al_i}$.

Given a map $f\colon X\to S$ and a polyhedron $Q\subset S$ a relative
polyhedron over $Q$ is a subset  $P$ with $\phi \colon U \to V \supset D_N$
such that $\phi\times f\colon U \to V\times S$ is a proper embedding
and $P= \phi^{-1}(D_N) \cap f^{-1}(Q)$.
A polyhedral covering $\cP$ of $X/S$ over a covering $\cQ$ of $S$ is a 
covering by relative polyhedra such that each $P_\al$ lies over some
$Q_{\nu(\al)}$. The mapping $\nu$ between index sets induces a
morphism $\nu \colon \cN(\cP)\to \cN(\cQ)$.

A {\sl resolving sheaf\/} for $X/S$ is a functor $\cR$ defined on 
$\cN(\cP)$ with values in the category of sheaves of graded 
differential algebras, such that each $\cR_A$ is a free graded
commutative algebra with a distinguished system of  generators
$e(A)$ of negative grading with differential $s$ such that the
complex $(\cR_A,s)$ is a resolution of the sheaf 
$$
\sier A:= (\phi_A\times f)_* \sier{}(X)|_{\wl D_A\times \wl Q_{\nu(A)}}\;.
$$
Furthermore, for $B\subset A$ each generator of $\cR_B$ is mapped
to a generator of $\cR_A$ under the morphism $(p_B^A)^*\cR_B\to
\cR_A$, where $p_B^A\colon \C^{N_A} \to \C^{N_B}$ is the canonical
projection, and the map $e(B)\to e(A)$ is injective.     
The elements in $e(A)$ which are not in the image of any of the
maps $e(B)\to e(A)$ with $B$ a proper subsimplex, are called
{\sl proper generators\/}. The functor $\cR$ is completely
determined by giving all the proper generators $e_{A,j}$ and the
values $s(e_{A,j})$.
The proof of the  existence of a resolving sheaf for $X/S$
for any polyhedral covering $\cP$ of $X$ over a covering $\cQ$ of 
$Y$ in \cite[Thm.~1.1]{P2} follows closely the absolute case
in \cite{P4}.

Let $\rho_A \colon D_A\times  Q_{\nu(A)}\to Y$. The functor $R$ from
$\cN(\cP)$ into the category of graded $\sier Y$ algebras, whose value
on $A$ is $R_A:= ({\rho_A})_*(\cR_A|_{D_A\times  Q_{\nu(A)}})$ is
called the {\sl resolution\/} of $X/S$ on $\cP$.

Let now $Y \mapright \pi X \mapright f S$ and $\cF$ a sheaf on $Y$.
Let $R$ be a resolution of $X/S$ on some polyhedral covering $\cP$,
$\cM$ a polyhedral covering of $\pi$ over $\cP$ and 
$\kappa\colon \cN(\cM)\to \cN(\cP)$ the induced map between the
nerves of the coverings. The functor $F=\pi_*\cF(\cM)$ has value
$\pi_*\cF(M_A)$ on $A\in \cN(\cM)$. Finally $f_*F$ is a functor from
$\cN(\cM)$ to the category of $\sier S$-sheaves.

\roep Definition{\/ \rm \cite[Def.~2.2]{P2}}.
The {\sl tangent complex\/} $T^*(R,F)=\sum_0^\infty T^n(R,F)$ of the
resolution $R$ with values in $\cF$ has as term of degree $n$ the 
$\sier S$-sheaf of $\sier S$-derivations of functors 
$\kappa^* R \to f_*F$ of degree $n$. In other words,
an element $v \in T^n(R,F)_s$ 
is a collection of compatible
$\sier{S,s}$-derivations $v_A\colon R_{\kappa(A),s}\to
(f\pi)_*(\cF|\cM_A)_{s}$ for $A\in \cN(\cM)$. 
The differential in the tangent complex is given by 
$dv=-(-1)^{\deg v}vs$ with $s$ the differential in $R$.
The cohomology of this complex is $T^n(X/S;\cF)$.
\endroep

\neu{spec-seq}
Let now $f\colon (X,x)\to (S,s)$ be  a finite map of germs.
We can choose the coverings 
$\cP$ of $X$ and $\cQ$ of $S$ to consist of one element each.

\proclaim Proposition.
There exists a spectral sequence
$$
E_2^{p,q}=T^p(X/S;R^q\pi_*\cF)
\Rightarrow T^n(X/S; \cF)\;.
$$

\roep Proof.
To compute the tangent complex $T^*(R,F)$ we describe the
resolving sheaf $\cR$ in more detail. Let $P_a$ be the polyhedron
covering $X$. Let $\cR_a$ be a resolution of $\sier a$ with
generators $e_{a,j}$. This is basically the Tyurina resolvent of $\sier {X,x}$.
The only simplices occurring in $\cN(\cP)$ are of the form
$A=(a,\dots,a)$. The proper generators of $\cR_A$ are
divided into a basic group, consisting of $e_{A,j}$ with
$\deg e_{A,j} = \deg e_{a,j} - \dim A$, and a complementary group
of elements of degree $-\dim A$,
corresponding to the coordinate functions $z_j$ on $\C^{n_a}$.
By considering the $z_j$ as complementary generators of $\cR_a$
of degree 0 we
have the relation $\deg e_{A,j} = \deg e_{a,j} - \dim A$
for all proper generators.
The improper generators of $\cR_A$ are of the form $e_{B,j}$ with
$B\subset A$ a proper subsimplex.

We introduce a filtration $\psi$ on $e(A)$ by setting 
$\psi(e_{B,j})=\deg e_{B,j} +\dim B$, which is the degree of the 
corresponding generator of $\cR_a$. This is the same filtration as in
the proof of \cite[Thm.~1.1]{P2}, but this special case is simpler because
there are no elements with filtration 1. We denote for all $k\leq0$
by $\cS_A^k$ the 
subalgebra of $\cR_A$ generated by the identity and the generators
of filtration at least $k$.
We define a derivation $\partial_A$ by
$$
\partial_A(e_{B,j})=(-1)^me_{B\setminus\beta_0,j}
+(-1)^{m-1}e_{B\setminus\beta_1,j}+\cdots+e_{B\setminus\beta_m,j},
$$
where $B=(\beta_0,\dots,\beta_m)$.
The proof of \cite[Thm.~1.1]{P2} yields that the differential 
in $\cR_A$ can be written as $s_A=\partial_A+g_A$ with $g_A$ a derivation
such that $g_A(e_{A,j})$ is a section of $\cS_A^{k+1}$ if 
$\psi(e_{A,j})=k$. 

Note that our $\partial_A$ differs from that
in \cite{P2}, in that it extends to the last index.
We indicate the first few terms of the resolution. 
Let $A=(a,a)$. To distinguish between
the indices we write $A=(\al,\beta)$. On $\C^{N_a}$ we have coordinates
$z_j$. The map $\phi_A$ embeds $U_A=U_a$ in the diagonal of
$\C^{N_a}\times \C^{N_a}$, so if we take coordinates $z_j^{(\al)}$
and $z_j^{(\beta)}$, we get the $N_a$  equations $z_j^{(\al)}
-z_j^{(\beta)}$, which equals $\partial_A(z^{(A)}_j)$. To each generator
$f$ of the ideal of $X$ correspond a generator $e_{A,f}$.
One has $s_A(e_{\al,f})=f(z^{(\al)})$, $s_A(e_{\beta,f})=f(z^{(\beta)})$
and $f(z^{(\beta)})-f(z^{(\al)})$ lies in the ideal of the diagonal,
so equals $\sum_j (z_j^{(\al)}-z_j^{(\beta)})h_j^{(A)}$
for some functions $h_j^{(A)}$.
Therefore $s_A(e_{A,f})= -e_{\beta,f}+e_{\al,f}+\sum_j h_j^{(A)} z^{(A)}_j$.

The filtration $\psi$ induces a decreasing
filtration on $T^*(R,F)$, of derivations
$\kappa^* R \to f_*F$ vanishing on the functor $\kappa^* S^{-k}$.
In its spectral sequence the differential $d_0$ is the differential
of the \v Cech complex of the covering $\cM$ of $Y$. The
$E_2$ term is therefore $T^p(X/S;R^q\pi_*\cF)$.
\qed

\section{The formula}
\neu{blow-exact}
We apply the results of the preceding section in the case that $X$ is a 
rational surface singularity of multiplicity $d$
and $\pi\colon \wtx\to X$ is the
first blow-up. Then $\pi_*\sier \wtx= \sier X$
and $R^q\pi_*(\sier \wtx)=0$ for $q>0$.
We choose a Noether normalisation 
$X\to S\cong \C^2$.
Note that the induced map  $\wtx\to \wt S$
to the blow-up of the plane is again a finite map of degree $d$,
so we can  compute the higher cotangent cohomology sheaves of $\wtx$ also
as relative Harrison cohomology:
$\cT^i_{\wtx} \cong \cT^i(\wtx/\wt S;\sier \wtx)$
for $i\geq2$.

By Prop.~\ref{spec-seq} we obtain $T^i(X/S;\sier \wtx)=T^i(X/S;\sier X)$,
which equals $T^i_X$ if $i\geq 2$.
The local-to-global spectral sequence  becomes
$$
E_2^{p,q}=H^p\left(\wtx,\cT^q(X/S; \sier \wtx)\right)
\Rightarrow T^n(X/S;\sier X)\;.
$$
The sheaves on  $\wtx$  occurring here are the basic objects of study.
To shorten notation we write
$$
\cF^i=\cT^{i}(X/S;\sier \wtx) \;.
$$
As only $H^0$ and $H^1$ contribute to the spectral sequence, we obtain
(cf.~\cite[Cor.~1.7]{CG}):

\proclaim Proposition.
If $\pi\colon \wtx\to X$ is the
first blow-up of a rational surface singularity
then one has short exact sequences
$$
0 \lra H^1(\wtx,\cF^{i-1}) \lra
T^i(X/S;\sier X) \lra H^0(\wtx,\cF^{i}) \lra 0 
$$
for all $i\geq0$.

\roep Remark.
{\sc Christophersen} and {\sc Gustavsen} \cite{CG} consider the absolute case. 
For $i\geq 2$ 
one has that $\cF^i\cong \cT^{i}(X;\sier \wtx)$.
Comparison of the exact sequences for $T^2_X$ and $T^2(X/S; \sier X)$
shows that even
$H^1(\wtx,\cF^1) \cong H^1(\wtx,\cT^1(X;\sier \wtx))$.
This can be seen directly as follows:
we have the exact sequence
$$
\cT^0(S;\sier \wtx)
\lra \cF^1 \lra \cT^1(X;\sier \wtx)\lra 0 \;.
$$
As $S$ is smooth 2-dimensional the sheaf
$\cT^0(S;\sier \wtx)$ is isomorphic to
$\sier \wtx \oplus \sier \wtx$ and therefore its $H^1$ vanishes,
implying the isomorphism in question.  

\neu{max-ideal}
Let $C$ be the exceptional curve with its scheme
structure, i.e., defined by $\m\sier \wtx$.
We study the sheaves  $\cF^i$  with the exact sequence
$$
0 \lra \m\cF^i \lra \cF^i \lra \cF^i_{|C} \lra 0 \;.
$$

\proclaim Proposition.
For all $i$ one has 
$$
\m \cF^i\cong \cT^i(\wtx/\wt S;\sier\wtx)\otimes \sier \wtx(iC) \;.
$$

\roep Proof.
We describe $X$ with coordinates $(x,y;z_1, \dots, z_m)$, where
$m=d-1$, and project onto the $(x,y)$-plane. The blow-up of the
plane can be covered by two charts. By genericity of the projection
we may assume that the two charts also suffice to cover $\wtx$. One chart
(which is a germ along the exceptional divisor, for
which we take a Stein representative) has coordinates
$(x,\eta;t_1, \dots, t_m)$ with $y=x\eta$, $z_i=xt_i$.
On the second chart we have  coordinates 
$(\xi,y;t_1', \dots, t_m')$ with $x=\xi y$, $z_i=yt_i'$
so on the intersection $t_i=\xi t_i'$.

We first look at the $x$-chart.
We compute the stalk $\cF^i_y$ from the Harrison complex
$C_{\rm sh}^\kb(\sier {X,0}/\sier {S,0};\sier{\wtx,y})$.
As abbreviation we write  $C_{\rm sh}^\kb(A/P;B)$.
An $(i+1)$-cochain 
$\vp$ is determined by its values
on $(i+1)$-tuples of coordinates $z_j$.
For $ \cT^i(\wtx/\wt S;\sier\wtx)$ we have the complex
$C_{\rm sh}^\kb(\sier{\wtx,y}/\sier{\wt S,\tilde{\pi}(y)};\sier{\wtx,y})$ 
(we write shortly $C_{\rm sh}^\kb(B/Q;B)$) 
with cochains determined by their values
on $(i+1)$-tuples of coordinates $t_j$.

We define maps $s\colon C_{\rm sh}^{i+1}(A/P,B) \to
C_{\rm sh}^{i+1}(B/Q,B)$ simply by
$$
(s\vp)(t_{j_0},\dots,t_{j_i})=\vp(z_{j_0},\dots,z_{j_i}) \;.
$$
This does not give a map of complexes, but the following
holds:
$$
sd=xds\;,
$$
where $x$ stands for the map `multiplication by $x$'.
We check:
$$
\def\qv{\qquad\qquad}
\displaylines{\qquad
(sd\vp)(t_{j_0},\dots,t_{j_i})=
(d\vp)(z_{j_0},\dots,z_{j_i})
\hfill\cr\qv{}=
z_{j_0} \vp (z_{j_1}, \dots,z_{j_i}) -
\vp (z_{j_0}z_{j_1}, z_{j_2},\dots,z_{j_i}) 
+\cdots\hfill\cr\qv{}=
z_{j_0} \vp (z_{j_1}, \dots,z_{j_i}) -
\sum p_{j_0j_1}^\nu(x,y)\vp (z_\nu, z_{j_2},\dots,z_{j_i}) 
+\cdots\hfill\cr\qv{}=
xt_{j_0} (s\vp) (t_{j_1}, \dots,t_{j_i}) -
\sum x{\wt p}_{j_0j_1}^\nu(x,\eta)(s\vp)  (t_\nu, t_{j_2},\dots,t_{j_i}) 
+\cdots\hfill\cr\qv{}=
xt_{j_0} (s\vp) (t_{j_1}, \dots,t_{j_i}) -
x(s\vp) (t_{j_0}t_{j_1}, t_{j_2},\dots,t_{j_i}) 
+\cdots\hfill\cr
\hfill{}=(xsd\vp)(t_{j_0},\dots,t_{j_i}) \;.\qquad\qv }
$$
The maps $s$ induce a map $s_*$ on cohomology:
if $d\vp=0$, then $x\,ds\vp=0$ and therefore
$ds\vp=0$ as $x$ is not a zero-divisor; furthermore, if
$\vp =d\psi$ then $s\vp=dxs\psi$.
As the maps $s$ are surjective, $s_*$ is also surjective.

The kernel of $s_*$ is the kernel of multiplication by $x$ on
$\Harr^{i+1}(A/P,B)$: if $s_*[\vp]=0$, then $s\vp=ds\psi$ for some $\psi$ and
therefore $x\vp=d\psi$ and conversely if $x\vp=d\psi$ then
$xs\vp=xds\psi$, so $s_*[\vp]=0$.
This makes $\Harr^{i+1}(B/Q,B)$ isomorphic to the image
$x\Harr^{i+1}(A/P,B)$, giving the claimed
isomorphism locally.

To see what happens globally we also look at the $y$-chart,
where we have the map $s'$.
If $\vp\in \m\cF^i$
we can write 
$\vp(z_{j_0},\dots,z_{j_i})=x\psi(z_{j_0},\dots,z_{j_i})$ 
and under the isomorphism it is mapped unto 
$(s\psi)(t_{j_0},\dots,t_{j_i})=\psi(z_{j_0},\dots,z_{j_i})$.
On the intersection we can write $\vp(z_{j_0},\dots,z_{j_i})
=y\xi\psi(z_{j_0},\dots,z_{j_i})$
so the other isomorphism maps to a cocycle homologous
to $(\xi s'\psi)(t_{j_0}',\dots,t_{j_i}')=\xi\psi(z_{j_0},\dots,z_{j_i})$. 
By $\xi$-linearity
we have that $(s\psi)(t_{j_0}',\dots,t_{j_i}')=
\xi^{i+1}(s\psi)(t_{j_0},\dots,t_{j_i})=
\xi^i (\xi s'\psi)(t_{j_0}',\dots,t_{j_i}')$. \qed

\roep Remark.
If $i\geq2$
then $\cT^i(\wtx/\wt S;\sier\wtx)\cong \cT^i_{\wtx}$ is concentrated in 
points, but in general $\cT^1(\wtx/\wt S;\sier\wtx)$ is not.
Its support is the critical locus of the map $\wtx \to \wt S$.
This follows from a local computation in smooth points of $\wtx$ 
using the exact sequence
$$
\cT^0_{\wtx}
\lra \cT^0(\wt S;\sier{\wtx}) \lra \cT^1(\wtx/\wt S;\sier\wtx)
\lra \cT^1_{\wtx}\lra 0\;. 
$$
In particular, if $C$ is reduced (i.e., $X$ has reduced
fundamental cycle), then the support has no one-dimensional
compact components and the  $H^1$ of the sheaf
vanishes.

%
\bigskip
\refstepcounter{subsection} \label{beperk}
\proclaim (\thesubsection)\enspace Proposition.
$
H^0(\wtx, \cF^i_{|C})=0
$.

\roep Proof. 
A global section of the sheaf $\cF^i_{|C}$ consists of a collection of
local sections $\vp \in  C_{\rm sh}^{i+1}(A/P,B)$ with $d\vp=0$
such that the difference between two of them lies in
the ideal of $C$. So each $\vp(z_{j_0},\dots,z_{j_i})$ gives rise
to a global section of $\sier C$, which therefore is a constant.
To prove the proposition we have to show that this constant
is zero.

We now argue as in the proof of Lemma  \ref{prel-zero}.
In the ring $B$ we have the equality
$$
\displaylines{
\qquad
0=(\delta f)(z_{j_0}, \dots, z_{j_\kn},z_k)=
z_{j_0} \vp (z_{j_1}, \dots,z_{j_i}) \pm z_kf(z_{j_0}, \dots, z_{j_\kn})
\hfill\cr\hfill{} 
-\sum p_{j_0j_1}^\nu(x,y)\vp (z_\nu, z_{j_2},\dots,z_{j_i})+\cdots 
\qquad\qquad
}
$$
and we may divide by $x$ to obtain
$$
\displaylines{
\qquad
0=
t_{j_0} \vp (z_{j_1}, \dots,z_{j_i}) \pm t_kf(z_{j_0}, \dots, z_{j_\kn})
\hfill\cr\hfill{} 
-\sum {\wt p}_{j_0j_1}^\nu(x,\eta)\vp (z_\nu, z_{j_2},\dots,z_{j_i})+\cdots
\;. 
\qquad\qquad
}
$$
We note that $C$ is a principal divisor, defined by $x$.
As $1$, $t_1$, \dots, $t_{d-1}$ are linearly independent
modulo $(x,\eta)$ we find that  $ \vp (z_{j_1}, \dots,z_{j_i})
=0 \in H^0(\sier C)$. 
\qed

\proclaim Corollary.
For $i\geq 2$ one has $H^0(\cF^i)\cong H^0(\cT^i_{\wtx})$.

\neu{kernel}
In order to compute $H^1(\cF^{i-1})$ we identify 
this group with the kernel of the map
$T^i(X/S;\sier X)\to H^0(\cF^i)$. From now on we fix an element $x$, which is 
supposed to be chosen generically.

\proclaim Lemma. The kernel of multiplication by $x$ on $T^i(X/S;\sier X)$ is
contained in the kernel of the map  $T^i(X/S;\sier X) \to H^0(\cF^i)$.

\roep Proof.
Suppose $\vp \in T^i(X/S;\sier X)$ is annihilated by $x$, i.e. 
$x \vp = d \psi$.
As $d\psi\equiv 0 $ modulo $P$, the argument of Lemma \ref{prel-zero}
shows that $\psi$ takes values in the maximal ideal.

Consider $\vp$ as global section of $\cF^i$. We may assume by
genericity of $x$ that a global section vanishes if and only if
it vanishes in the chart $x\neq0$ (the special points of $\cF^i$
on the exceptional divisor lie in this chart). 
The values of $\psi$ lie in $\m$ and are therefore in $B$ divisible
by $x$. We obtain that $\vp =d (\psi/x)$.
\qed

\proclaim Lemma. Let $K$ be a submodule of $T^i_X$, $i\geq 2$ 
containing the kernel of multiplication by $x$. Then $\dim K/xK=\cg T^i$.

\roep Proof.
Consider the multiplication $K\mapright{\cdot x} K$. As $K\subset T^i$
the kernel is always contained in $\Ker\{T^i\mapright{\cdot x} T^i\}$
so under the assumption of the Lemma both kernels are equal.
As $K$ is finite-dimensional,  kernel and  cokernel have the same
dimension, and the same holds for $T^i$. Therefore 
$\dim K/xK=\dim T^i/xT^i$. For a general hyperplane section
this dimension is  $\cg T^i$.
\qed

The previous two lemmas show that $\dim H^1(\cF^{i-1})/
xH^1(\cF^{i-1})$ has dimension $f_i(d)$. We have to determine
$xH^1(\cF^{i-1})$.
Since $x$ is generic the cokernel of
$\cF^{i-1}\mapright{\cdot x}\m \cF^{i-1}$ has support at the strict
transform of the divisor of $x$. In particular its $H^1$ vanishes
and therefore
the map $H^1(\cF^{i-1}) \mapright{\cdot x} H^1(\m \cF^{i-1})$
is surjective.
The image of the composed map $H^1(\cF^{i-1})\mapright{\cdot x}
H^1(\m \cF^{i-1}) \to H^1(\cF^{i-1})$ is $xH^1(\cF^{i-1})$.
As $H^1(\m \cF^{i-1}) \to H^1(\cF^{i-1})$ is 
injective we obtain $xH^1(\cF^{i-1})= H^1(\m \cF^{i-1})$.
For  $i> 2$ 
this group vanishes.

\neu{result}
To collect our results in formulas for the dimension of $T^i$
we introduce a name for the remaining unknown term.

\roep Definition.
For a rational surface singularity $X$ we define the
invariant
$$
c(X) = \dim H^1(\m \cF^1) \;.
$$
\endroep

\roep Remarks.\\
1)
By Remark \ref{max-ideal} the correction term $H^1(\m \cF^1)= 
H^1(\wtx,\cT^1(\wtx/\wt S;\sier\wtx)(C))$
vanishes if  $X$ has reduced fundamental cycle.
\\
2) 
Our invariant is the  same as the one {\sc Christophersen--Gustavsen} 
\cite{CG} define using their absolute version $\cF^1_{CG}$ 
of the sheaf $\cF^1$.
This follows from the formula for $T^2$. A direct proof can
be obtained as in Remark \ref{blow-exact}. One  tensors  
the exact sequence occurring there
with the invertible sheaf $\sier \wtx(-C)$
and uses  that $H^1(\wtx, \m\cF^1_{CG})=
H^1(\wtx, \cF^1_{CG}(-C))$ \cite[Prop.~4.1]{CG}.

\proclaim Theorem.
If $X$ is  a rational surface singularity of multiplicity $d$
and $\wtx\to X$ the first blow up, then for $i>2$
$$
\dim T^i_X= \cg T^i_X + \dim T^i_{\wtx}
$$
with $\cg T^i_X=f_i(d)$ as given in\/ \mref{prel-rnc}, and
$$
\dim T^2_X= (d-1)(d-3) + \dim T^2_{\wtx} +c(X) \;.
$$


\section{Good maximal deformations}

\neu{codim}
The dimension formula above is an inductive formula. To make it more 
explicit we first define
the {\sl multiplicity sequence\/} of a rational singularity
in the obvious way as the sequence of multiplicities of
the singularities on successive blow ups 
(this are the infinitely near singularities,
including the singularity itself). We denote by $X_P$ the singularity
at an infinitely near point $P$ and by $d(P)$ its multiplicity.

The $T^1$-formula of \cite{CG} can be best stated as a formula for
the codimension of the (smooth) Artin component in the Zariski tangent space
of the versal base space.
We denote this by invariant by $\cod(X)$. For the cone over the rational
normal curve of degree $d\geq3$ it has value $d-3$.

By induction we obtain from \cite[Thms.~3.8 and 3.13]{CG} with 
Thm.~\ref{result}  
the following formulas.

\proclaim Theorem.
For a rational surface singularity $X$ of multiplicity $d$
\begin{eqnarray*}
\dim T^i_X&=&\sum_P f_i(d)\;,\qquad i\geq 3\;,
\\
\dim T^2_X&=&\sum_P (d(P)-1)(d(P)-3) + \sum_P c(X_P)\;,
\\
\cod(X)&=& \sum_P(d(P)-3) + \sum_Pc(X_P) \;,
\end{eqnarray*}

where the sum ranges over all infinitely near singular points $P$
of multiplicity at least three.

\neu{dJvS}
{\sc De Jong} and {\sc Van Straten} \cite{dJvS} derived 
their dimension formulas for $T^1$ and $T^2$ 
using a special deformation  to the cone over the rational normal curve
of degree $d$ and all singularities on the first blow up. 
The same deformation also yields the surjectivity
of the obstruction map.

The important ingredient in our  dimension formulas is the
multiplicity sequence. With this in mind we define a more general class
of special deformations.

\roep Definition.
A {\sl good maximal deformation\/} of a rational surface singularity
$X$ of multiplicity $d$ 
is a one-parameter deformation $X_T\to T$ such that the
general fibre $X_t$ has as singularities
cones over rational normal curves of multiplicity $d(P)$, 
one for each infinitely near singularity of multiplicity at least 3.
\endroep

\proclaim Proposition.
If a rational surface singularity $X$ has a good maximal deformation,
then the number of equations of the versal base space equals
the dimension of $T^2_X$.

\roep Proof.
Let $X_T \to T$ be a good maximal deformation  of $X$.
Versality induces a map $b\colon T \to B_X$,  to the base space of the
miniversal deformation.
In a general point $b(t)$ 
the Zariski tangent space to $B_X$ has codimension 
$\sum_P c(X_P)$ in the Zariski tangent space at the origin, so one needs
so many linear equations. To describe the base space at the point $b(t)$
one needs in addition $\sum_P (d(P)-1)(d(P)-3)$ quadratic
equations. Therefore one needs at the origin  at least 
$\sum_P (d(P)-1)(d(P)-3) + \sum_P c(X_P)$ equations.
As this number is the dimension of $T^2_X$, they also suffice.
\qed

\roep Remark.
The $T^1$-formula above follows easily from
the existence of a good maximal deformation.
Consider the long exact sequence 
$$
T^1_{X_T/T} \mapright{\alpha} T^1_X
\lra
T^2_{X_T/T} \mapright{\cdot t} T^2_{X_T/T}  \mapright{\beta} T^2_X
\lra
T^3_{X_T/T} \mapright{\cdot t} T^3_{X_T/T} \mapright{\gamma} T^3_X
$$
The dimension
of $T^3$ is constant in this deformation,
so the rank of the $\C\{t\}$-module $T^3_{X_T/T}$ is 
equal to $\dim T^3_X$. Therefore $\gamma$ is surjective,
and as the rank is equal to  $\dim \Coker(\cdot t) - \dim \Ker(\cdot t)$, we
obtain that multiplication by $t$ is injective, and therefore
the map $\beta$ is also surjective.
By assumption on the deformation we have
$\dim T^2_X-\Rank T^2_{X_T/T}=\sum_P c(X_P)$.
The dimension of $\Ker(\cdot t\colon T^2_{X_T/T} \lra T^2_{X_T/T})$
equals therefore also $\sum_P c(X_P)$ and
we obtain $\dim T^1_X= \dim {\rm Im}\; \alpha + \sum_P c(X_P)$.
But we  know the dimension of the image of $\alpha$: by \cite{GL}
it equals the dimension of the Zariski tangent space 
to $B_X$  in the general point $b(t)$. The codimension
of the Artin component is there $\sum_P (d(P)-3)$.
\endroep

\neu{niet}
The existence of a good maximal deformation is established
by {\sc De Jong} and {\sc Van Straten}  for rational
singularities  with reduced
fundamental cycle and by {\sc De Jong} in the determinantal case
\cite{dJ}.

In the first case one can deform to the singularities  on
the first blow up plus a cone of degree $d$, but for determinantal
singularities this is in general impossible.
Specifically, if $X$ contains one configuration
of type $D_{2d+1}^{\rm II}$ (in the notation of \cite{dJ}), then
the dimension of $T^1$ will be too small to allow a deformation
to a determinantal with $A^1_{2d}$-configuration and a cone
(I checked this by computing $T^1$ in the case $d=3$). 
But the  singularity does deform into two
cones.

In the definition of a good maximal deformation we ignore all occurring
rational double points, for a good reason:
the singularities on the first blow-up
of $D_4$
are three $A_1$'s, and the cone
over the rational normal curve of degree two is also $A_1$
but there
is no deformation  $D_4 \to 4A_1$.

I hoped to prove that every   rational singularity has a good
maximal deformation but unfortunately this is not true.
First we note that the dimension 
of the Artin component equals
$h^1(\widetilde X, \Theta_{\widetilde X})$, where 
$\widetilde X \to X$ is the minimal resolution, with exceptional set
$E=E_1\cup \dots \cup E_r$. As usual we denote by $-b_i$ the 
self-intersection of the irreducible component $E_i$.
By 
\cite[Prop.~2.2 and 2.5]{W}, 
$$
 h^1(\widetilde{X}, \Theta_{\widetilde X}) = \sum_{i=1}^r(b_i-1) +
 h^1(\widetilde{X}, \Theta_{\widetilde X}(\log E)).
$$
The second summand gives the dimension of the equisingular stratum.
The stratum in the Artin component of fibres with a cone over
the rational normal curve of multiplicity $d$ as singularity
has codimension $d-1$. By openess of versality these strata intersect
transversally in the base of a good maximal deformation.

\proclaim Lemma. 
The general rational singularity with a given resolution graph
does not have a good maximal deformation, if\/
$
\sum_P (d(P)-1) \geq \sum (b_i-1)
$.

An example where this condition is satisfied, 
is obtained by generalising the
$D_4$-singularity to higher multiplicity:
consider a singularity with fundamental cycle reduced everywhere
except at one $(-2)$-curve, such that the first blow up has
three singularities of multiplicity $k$.
The simplest way to do this gives the following graph,
where as usual a dot stands for a $(-2)$-curve.
$$
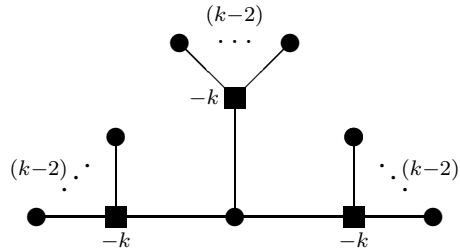
\def\dstod{\mathinner{\mkern1mu\raise1pt\hbox {.}\mkern 2mu 
 \raise 5.5pt\hbox {.}\mkern 2mu\raise10pt\vbox{\kern10pt\hbox{.}}\mkern 1mu}}
\def\ddots{\mathinner{\mkern1mu\raise10pt\vbox{\kern10pt\hbox{.}}\mkern 2mu 
 \raise 5.5pt\hbox {.}\mkern 2mu\raise1pt\hbox {.}\mkern 1mu}}
\unitlength=30pt
\def\ci{\circle*{0.23}}
\def\vi{{\makebox(0,0){\rule{8pt}{8pt}}}}
\def\mbt#1{\makebox(0,0)[t]{$\scriptstyle #1$}}
\begin{picture}(5,3.5)(-2.5,-.5)
\put(-2.5,0){\line(1,0){5}}
\put(-2.5,0){\ci}
\put(-1.5,0){\vi}
\put(-1.5,-.2){\mbt{-k}}
\put(1.5,0){\vi}
\put(1.5,-.2){\mbt{-k}}
\put(2,.65){\makebox(0,0){$\ddots$}}
\put(2.1,.6){\makebox(0,0)[l]{$\scriptstyle (k-2)$}}
\put(-2,.65){\makebox(0,0){$\dstod$}}
\put(-2.1,.6){\makebox(0,0)[r]{$\scriptstyle (k-2)$}}
\put(1.5,0){\line(0,1){1}}
\put(-1.5,0){\line(0,1){1}}
\put(-1.5,1){\ci}
\put(1.5,1){\ci}
\put(0,0){\line(0,1){1.5}}
\put(0,1.5){\vi}
\put(-.2,1.5){\makebox(0,0)[r]{$\scriptstyle -k$}}
\put(0,1.5){\line(1,1){.7}}
\put(0,1.5){\line(-1,1){.7}}
\put(.7,2.2){\ci}
\put(-.7,2.2){\ci}
\put(0,2.4){\makebox(0,0)[b]{$\scriptstyle (k-2)$}}
\put(0,2.2){\makebox(0,0)[b]{$\ldots$}}
\put(1.5,1){\ci}
\put(0,0){\ci}
\put(2.5,0){\ci}
\end{picture}
$$
The singularity has multiplicity $3k-4$;
for $k=2$ it is indeed $D_4$.
Here we have $\sum_P (d(P)-1) = \sum (b_i-1)=6k-8$.

\parskip=0pt plus 1pt

\vfill

{\obeylines
Address of the author:
Matematik
G\"oteborgs universitet
Chalmers tekniska h\"ogskola
SE 412 96 G\"oteborg, Sweden
email: stevens@math.chalmers.se}

\end{document}